\documentclass[11pt]{article}
\title{ A Note on M\"{o}bius Function and M\"{o}bius Inversion Formula of Fibonacci Cobweb Poset}
\author{Ewa Krot \\
\\Institute of Computer Science, Bia{\l}ystok University\\
PL-15-887 Bia{\l}ystok, ul.Sosnowa 64, POLAND\\
e-mail: ewakrot@wp.pl, ewakrot@ii.uwb.edu.pl}

\usepackage{amsmath,amsthm}

\chardef\bslash=`\\ 
\newtheorem{thm}{Theorem}[section]

\begin{document}
\maketitle
\begin{abstract}
The explicit formula for M\"{o}bius function of Fibonacci cobweb
poset $P$ is given here for the first time by the use of
Kwa\'sniwski's definition of $P$ in plane grid coordinate system
\cite{1}.
\end{abstract}

\section{Fibonacci cobweb poset}
The Fibonacci cobweb poset $P$ has been introduced by
A.K.Kwa\'sniewski in \cite{3,4} for the purpose of finding
combinatorial interpretation of fibonomial coefficients and their
reccurence relation. At first the partially ordered set $P$
(Fibonacci cobweb poset) was defined
 via its Hasse diagram as follows: It looks
like famous rabbits grown tree but it has the specific cobweb in
addition, i.e. it consists of levels labeled by Fibonacci numbers
(the $n$-th level consist of $F_{n}$ elements). Every element of
$n$-th level ($n \geq 1, n \in {\bf N}$) is in partial order
relation with every element of the $(n+1)$-th level but it's not
with any element from the level in which he lies ($n$-th level)
except from it.

 \section{The Incidence Algebra ${\bf I(P)}$}
One can define the incidence algebra of ${\bf P}$ (locally finite
partially ordered set) as follows (see \cite{7,8}):
$$ {\bf I(P)}=\{f:{\bf P}\times {\bf P}\longrightarrow {\bf R};\;\;\;\; f(x,y)=0\;\;\;unless\;\;\; x\leq y\}.$$
The sum of two such functions $f$ and $g$ and multiplication by
scalar are defined as usual. The product $H=f\ast g$ is defined as
follows:
$$ h(x,y)=(f\ast g)(x,y)=\sum_{z\in {\bf P}:\;x\leq z\leq y} f(x,z)\cdot g(z,y).$$
It is immediately verified that this is an associative algebra
over the real field ( associative ring).

The incidence algebra has an identity element $\delta (x,y)$, the
Kronecker delta. Also the zeta function of ${\bf P}$ defined for
any poset by:
$$ \zeta (x,y)=\Big\{\begin{array}{l}1\;\;for\;\;\; x \leq y\\0\;\;otherwise\end{array}$$
 is an element of ${\bf I(P)}$. The one for Fibonacci cobweb poset  was expressed by $\delta$ in
\cite{1,4} from where quote the result:
\begin{equation}\label{dzeta}
\zeta =\zeta_{1}-\zeta_{0}
\end{equation}
 where for $x,y \in {\bf N}$:
\begin{equation}
\zeta_{1}(x,y)=\sum_{k=0}^{\infty}\delta(x+k,y)
\end{equation}
\begin{equation}
\zeta_{0}(x,y)=\sum_{k \geq 0}\sum_{s \geq 0}\delta
(x,F_{s+1}+k)\sum_{r=1}^{F_{s}-k-1}\delta (k+F_{s+1}+r,y).
\end{equation}

The knowledge of $\zeta$ enables us to construct other typical
elements of incidence algebra perfectly suitable for calculating
number of chains, of maximal chains etc. in finite sub-posets of
${\bf P}$. The one of them is M\"{o}bius function indispensable in
numerous inversion type formulas of countless applications. It is
known that the $\zeta$ function of a locally finite partially
ordered set is invertible in incidence algebra and its inversion
is the famous M\"{o}bius function $\mu$ i.e.:
$$\zeta \ast \mu=\mu \ast \zeta=\delta.$$The M\"{o}bius function
$\mu$ of Fibonacci cobweb poset ${\bf P}$ was presented  for the
first time by the present author in \cite{2}. It was recovered by
the use of the recurrence formula for M\"{o}bius function of
locally finite partially ordered set  ${\bf I(P)}$ (see \cite{7}):
\begin{equation}
 \left\{\begin{array}{l}\mu(x,x)=1\;\;\;\;for\;\; all \;\;x\in{\bf
 P}\quad \quad \quad \\\ \\
\mu (x,y)=-\sum_{x\leq z<y}\mu (x,z)\end{array}\right.
\end{equation}

The  M\"{o}bius function of Fibonacci cobweb poset was there given
by following formula:
\begin{equation}\label{mobius1}
\mu(x,y)=\left\{\begin{array}{l} 0\;\;\;\;\;\;x>y\\
1\;\;\;\;\;\;x=y\\
0\;\;\;\;\;\;\;F_{k+1}\leq x,y\leq F_{k+2}-1;\;x\neq y;\;k\geq 3\\
\\-1\;\;\;\;\;\;F_{k+1}\leq x\leq F_{k+2}-1<F_{k+2}\leq y\leq
F_{k+3}-1\\ \\
(-1)^{n-k}\prod_{l=k+1}^{n-1}(F_{l}-1)\quad \quad \quad F_{k+1}\leq x\leq F_{k+2}-1,\\
\quad \quad \quad \quad \quad \quad \quad \quad \quad \quad
F_{n+1}\leq y\leq F_{n+2}-1;\;n>k+1,\end{array}\right.
\end{equation}
where:
\begin{itemize}
\item the condition $F_{k+1}\leq x,y\leq F_{k+2}-1;\;x\neq
y;\;k\geq 3$ means that $x,y$ are different elements of $k$-th
level; \item the condition $F_{k+1}\leq x\leq
F_{k+2}-1<F_{k+2}\leq y\leq F_{k+3}-1$ means that $x$ is an
element of $k$-th level and $y$ is an element of $(k+1)$-th level;
\item the condition $F_{k+1}\leq x\leq F_{k+2}-1,\;F_{n+1}\leq
y\leq F_{n+2}-1;\;n>k+1$ means that $x$ is an element of $k$-th
level and $y$ is an element of $n$-th level.
\end{itemize}

The above formula allows us to find out the $\mu$ function matrix
(see\cite{2})but it is not good enough to be applied in compact
form via  M\"{o}bius inversion formula for cobweb poset. For this
purpose more convenient, explicit formula  is needed.
\section{Plane grid coordinate system of $P$}
In \cite{1} A. K. Kwa\'sniewski defined cobweb poset $P$ as
infinite labeled graph oriented upwards as follows: Let us label
vertices of $P$ by pairs of coordinates: $\langle i,j \rangle \in
{\bf N}\times {\bf N}$, where the second coordinate is the number
of level in which the element of $P$ lies (here it is the $j$-th
level) and the first one is the number of this element in his
level (from left to the right), here $i$  . We shall refer,
(following \cite{1}) $\Phi_{s}$ as to the set of vertices
(elements) of the $s$-th level, i.e.:
$$\Phi_{s}=\left\{\langle j,s \rangle ,\;\;1\leq j \leq F_{s}\right\},\;\;\;s\in{\bf N}.$$
For example $\Phi_{1}=\{\langle 1,1\rangle\},\; \Phi_{2}=\{\langle
1,2\rangle\}, \;\Phi_{3}=\{\langle 1,3\rangle, \langle
2,3\rangle\},\\\Phi_{4}=\{\langle 1,4\rangle, \langle 2,4\rangle,
\langle 3,4\rangle\},\;\Phi_{5}=\{\langle 1,5\rangle, \langle
2,5\rangle, \langle 3,5\rangle\,\langle 4,5\rangle,\langle
5,5\rangle\}$ ....

Then $P$ is a labeled graph $P=\left(V,E\right)$ where
$$V=\bigcup_{p\geq1}\Phi_{p},\;\;\;E=\left\{\langle \,\langle j,p\rangle,\langle q,p+1
\rangle\,\rangle\right\},\;\;1\leq j\leq F_{p},\;\;1\leq q\leq
F_{p+1}.$$

Now we can define the partial order relation on $P$ as follows:
let\\ $x=\langle s,t\rangle, y=\langle u,v\rangle $ be elements of
cobweb poset $P$. Then
$$ ( x \leq y) \Longleftrightarrow
 [(t<v)\vee (t=v \wedge s=u)].$$

\section{The M\"{o}bius function and M\"{o}bius inversion formula
on $P$} The above definition of $P$ allows us to derive an
explicit formula for M\"{o}bius function of cobweb poset $P$. To
do this we can use the formula (\ref{mobius1}). Then for
$x=\langle s,t\rangle,\;\;y=\langle u,v\rangle,\,\,1\leq s\leq
F_{t},\;\;1\leq u\leq F_{v},\;t,v\in {\bf N}$ we have
\begin{multline}\label{mobius2}
\mu(x,y)=\mu\left( \langle s,t\rangle ,\langle u,v\rangle
\right)=\\
\qquad \quad \;
=\delta(t,v)\delta(s,u)-\delta(t+1,v)+\sum_{k=2}^{\infty}\delta(t+k,v)(-1)^{k}
\prod_{i=t+1}^{v-1}(F_{i}-1)
\end{multline}
where $\delta$ is the Kronecker delta defined by
$$ \delta (x,y)=\Big\{\begin{array}{l}1\;\;\;\;\; x=y\\0\;\;\;\;\;x\neq y\end{array}.$$
We can also derive more convenient then (\ref{dzeta}) formula for
$\zeta$ function of $P$ \\( the characterictic function of partial
order relation in $P$):
\begin{equation}
\zeta(x,y)=\zeta\left( \langle s,t\rangle ,\langle u,v \rangle
\right)=\delta (s,u)\delta (t,v)+\sum_{k=1}^{\infty}\delta(t+k,v).
\end{equation}
The formula (\ref{mobius2}) enables us to formulate following
theorem (see \cite{7}):
\begin{thm}{\bf (M\"{o}bius Inversion Formula of cobweb $P$)}\\
Let $f(x)=f(\langle s,t\rangle)$ be a real valued function,
defined for $x=\langle s,t\rangle$ ranging in cobweb poset $P$.
Let an element $p=\langle p_{1},p_{2}\rangle$ exist with the
property that $f(x)=0$ unless $x\geq p$.

Suppose that $$g(x)=\sum_{\left\{ y \in P :\; y\leq x
\right\}}f(y).$$ Then
$$f(x)=\sum_{\{y\in P:\;y\leq x\}}g(y)\mu(y,x).$$
But using coordinates of $x,y$ in $P$ i.e. $x=\langle
s,t\rangle,\;y=\langle u,v\rangle$ if
$$g(\langle s,t\rangle)=\sum_{v=1}^{t-1}\sum_{u=1}^{F_{v}}\left(f(\langle u,v\rangle )\right)+f(\langle s,t\rangle)$$
then we have
\begin{multline}
f(\langle s,t\rangle)=\sum_{v\geq 1}\sum_{u=1}^{F_{v}}g(\langle
u,v\rangle)\mu(\langle s,t\rangle ,\langle
u,v\rangle)=\\
=\sum_{v\geq 1}\sum_{u=1}^{F_{v}}g(\langle
u,v\rangle)\left[\delta(v,t)\delta(u,s)-\delta(v+1,t)+\sum_{k=2}^{\infty}\delta(v+k,t)(-1)^{k}
\prod_{i=v+1}^{t-1}(F_{i}-1)\right].
\end{multline}
\end{thm}

 AMS Classification numbers: 11C08, 11B37, 47B47
 \end{document}